\documentclass[twoside,a4paper,12pt,reqno]{amsart}
\usepackage{amssymb}
\usepackage{fullpage}
\usepackage{lmodern}

\newtheorem{theorem}{Theorem}[section]
\newtheorem{lemma}[theorem]{Lemma}
\newtheorem{proposition}[theorem]{Proposition}
\newtheorem{corollary}[theorem]{Corollary}

\theoremstyle{remark}

\numberwithin{equation}{section}
\usepackage[colorlinks,hypertexnames=false,linkcolor=blue,urlcolor=blue,citecolor=blue]{hyperref}
\hypersetup{pdftitle={Harish-Chandra isomorphism over arbitrary base rings},pdfauthor={Yu Shi},pdfsubject={A proof by the parabolic intersection lemma}}
\newcommand{\Z}{\mathbb Z}
\newcommand{\Q}{\mathbb Q}
\newcommand{\F}{\mathbb F}
\newcommand{\Gm}{\mathbb G_m}
\newcommand{\cO}{\mathcal O}
\newcommand{\fg}{\mathfrak g}
\newcommand{\fh}{\mathfrak h}
\newcommand{\fn}{\mathfrak n}
\newcommand{\fl}{\mathfrak l}
\newcommand{\fq}{\mathfrak q}
\newcommand{\fu}{\mathfrak u}
\newcommand{\fB}{\mathfrak b}
\newcommand{\cA}{\mathcal A}
\DeclareMathOperator{\Sym}{Sym}
\DeclareMathOperator{\Lie}{Lie}
\DeclareMathOperator{\Ad}{Ad}
\DeclareMathOperator{\HC}{HC}
\DeclareMathOperator{\gr}{gr}
\DeclareMathOperator{\Spec}{Spec}
\DeclareMathOperator{\Ann}{Ann}
\DeclareMathOperator{\Frac}{Frac}

\newcommand{\iso}{\xrightarrow{\ \sim\ }}
\newcommand{\into}{\hookrightarrow}

\title[Harish-Chandra isomorphism over arbitrary base rings]{Harish-Chandra isomorphism\\
over arbitrary base rings}

\author{Yu Shi}
\address{Yu Shi, Humboldt-Universit\"at zu Berlin, Institut f\"ur Mathematik, Unter den Linden 6, Berlin, 10099 (Germany)}
\email{yu.shi@student.hu-berlin.de}

\begin{document}

\begin{abstract}
Let \(G\) be a split connected reductive group scheme over \(\Z\), and let \(R\) be any commutative ring. We prove that Harish-Chandra projection induces a natural, strict filtered isomorphism
\[
U_R(\fg_R)^{G_R}\ \simeq\ \Sym_R(\fh_R)^{W,\bullet}
\]
whereas the corresponding Chevalley restriction need not be surjective even over a field. 

\end{abstract}

\maketitle

\tableofcontents


\section{Introduction}\label{sec:introduction}

\subsection{Main results}

The classical Harish-Chandra isomorphism identifies the center of the enveloping algebra of a split reductive Lie algebra over a characteristic-zero field with the invariants of the shifted Weyl action. We prove a strict filtered version for group-scheme invariants over every commutative coefficient ring and for every split connected reductive group scheme over \(\Z\).

Let \(G\) be a split connected reductive group scheme over \(\Z\), and fix a split maximal torus \(T\subset B\) in a Borel subgroup. Write
\[
\fg=\Lie(G),\qquad \fh=\Lie(T),\qquad W=N_G(T)/T.
\]
For a commutative unital ring \(R\), set
\[
G_R=G\times_{\Spec\Z}\Spec R,\qquad
U_R=U_R(\fg_R),\qquad S_R=\Sym_R(\fh_R).
\]
We equip the usual enveloping algebra \(U_R\) with its PBW filtration \(F_\bullet U_R\), and \(S_R\) with the polynomial-degree filtration. Superscripts \(G_R\) denote group-scheme invariants, as defined in Section~\ref{sec:preliminaries}.

Let \(\rho\) be the half-sum of the positive roots determined by \(B\). The dot action (or shifted action) of \(W\) on \(\fh^*\) is
\[
w\bullet\lambda=w(\lambda+\rho)-\rho.
\]
This action is integral because \(w\rho-\rho\) belongs to the root lattice, even when \(\rho\) is not a character of \(T\). It induces an action on \(S_R\). The triangular decomposition relative to \(B\) defines the PBW Harish-Chandra projection \(\HC_R:U_R\to S_R\).

\begin{theorem}\label{thm:main}
For every commutative unital ring \(R\), Harish-Chandra projection induces a natural isomorphism of \(R\)-algebras
\[
\HC_R:U_R^{G_R}\iso S_R^{W,\bullet}.
\]
It is strict for the induced filtrations: for every \(d\ge0\),
\[
\HC_R\bigl(F_dU_R\cap U_R^{G_R}\bigr)
=F_dS_R\cap S_R^{W,\bullet}.
\]
\end{theorem}

Throughout the paper, we assume that \(G\) is not a torus. In the torus case, \(U_R=S_R\), the adjoint and Weyl actions are trivial, and Harish-Chandra projection and Chevalley restriction are identity maps.
The following corollary brings the Harish-Chandra isomorphism into its more familiar form.

\begin{corollary}\label{cor:fullcenter}
For every commutative unital ring \(R\),
\[
Z(U_R)=U_R^{G_R}
\quad\Longleftrightarrow\quad
R\text{ is torsion-free over }\Z.
\]
Under these equivalent conditions, Harish-Chandra projection induces an isomorphism
\[
\HC_R:Z(U_R)\iso S_R^{W,\bullet}.
\]
The isomorphism is strict filtered and natural among torsion-free coefficient rings.
\end{corollary}
\begin{proof}
The equivalence follows from Lemma~\ref{lem:center}, and the isomorphism and its properties follow from Theorem~\ref{thm:main}.
\end{proof}

\subsection{Base change}\label{subsec:sl2example}

Naturality in Theorem~\ref{thm:main} means that Harish-Chandra projection commutes with changes
of coefficients \(R\to R'\). However,
the base change maps
\[
U_R^{G_R}\otimes_RR'\longrightarrow U_{R'}^{G_{R'}},\qquad
S_R^{W,\bullet}\otimes_RR'\longrightarrow S_{R'}^{W,\bullet}
\]
may not be isomorphisms. They are isomorphisms when \(R'\) is flat
over \(R\), since flat tensor product preserves the kernels defining
invariants. In particular, Corollary~\ref{cor:fullcenter} gives
\(Z(U_{\Z})\otimes_{\Z}R\simeq Z(U_R)\) for every torsion-free
\(\Z\)-algebra \(R\).

For example, take \(G=\mathrm{SL}_2\), with standard generators
\(e,f,h\), and put
\[
y=h+1,\qquad C=y^2+4fe,\qquad I_R=U_R^{G_R}.
\]
The dot reflection sends \(y\) to \(-y\).
For \(k=\F_2\) and \(A=\Z/4\Z\), the rank-one calculation in
Proposition~\ref{prop:rankone} gives
\begin{equation}\label{eq:intro-sl2-z4}
I_{\Z}=\Z[C],\qquad I_k=k[y],\qquad
I_A\simeq A[c,t]/(2t,t^2),
\end{equation}
where, over \(A\), \(c\) maps to \(C=y^2\) and \(t\) to \(2y\).
Reduction of integral invariants to \(k\) has image
\(k[y^2]\subsetneq k[y]\). Moreover, the map
\[
I_A\otimes_A k\simeq k[c,t]/(t^2)\longrightarrow I_k=k[y],
\qquad c\longmapsto y^2,\quad t\longmapsto0
\]
has kernel \((t)\) and image \(k[y^2]\), so it is neither injective nor
surjective. Thus an invariant can disappear under base change, while a
new one appears in the target.

\subsection{Chevalley restriction}\label{subsec:chevalley}

The PBW filtration leads to the (coadjoint) Chevalley restriction map
\begin{equation}\label{eq:intro-chevalley}
\operatorname{res}_{\fh}:
\Sym_R(\fg_R)^{G_R}\longrightarrow S_R^W.
\end{equation}
Here \(\Sym_R(\fg_R)=\cO(\fg_R^*)\), and \(W\) acts without the
shift. This map need not be surjective, even over a field.

Take \(G=\mathrm{PGL}_2\) and \(k=\F_2\). Let \(e,f\) be the
off-diagonal matrix units and \(H\) the class of
\(\operatorname{diag}(1,0)\). The ordinary reflection fixes \(H\).
With \(Q=ef+H^2\), the Chevalley restriction is
\begin{equation}\label{eq:intro-chevalley-counterexample}
k[Q]\longrightarrow k[H],\qquad Q\longmapsto H^2,
\end{equation}
whose image is \(k[H^2]\subsetneq k[H]\).
This failure is explicitly remarked in \cite[Rem.~3.11]{Spice}.
It also corresponds to the \(\mathrm{SL}_2\) example of adjoint Chevalley restriction in
\cite[Fact~6.1 and Prop.~6.2]{CR}: the trace pairing identifies
\(\mathfrak{pgl}_2^*\) with \(\mathfrak{sl}_2\), on which the adjoint
action of \(\mathrm{SL}_2\) factors through \(\mathrm{PGL}_2\).

Harish-Chandra projection nevertheless remains an isomorphism.
Indeed, the dot reflection sends \(H\) to \(H+1\), and
Proposition~\ref{prop:rankone} gives
\[
\HC_k:k[fe+H(H+1)]\iso k[H(H+1)].
\] 
Since \(H(H+1)\) has leading term \(H^2\), the degree filtration gives
\[
\gr\bigl(S_k^{W,\bullet}\bigr)=k[H^2]\subsetneq k[H]=S_k^W.
\]

A coroot is \emph{primitive} if it belongs to a \(\Z\)-basis of
\(X_*(T)\), equivalently if its image in
\(X_*(T)\otimes_\Z\F_p\) is nonzero for every prime \(p\).
For coroots, nonprimitivity is equivalent to divisibility by \(2\),
since the pairing with the corresponding root is \(2\).
Write \(R[2]=\Ann_R(2)=\{a\in R:2a=0\}\). Then there is a precise criterion:

\begin{theorem}

    \(
\operatorname{res}_{\fh}\) is an isomorphism if and only if 
\(
R[2]=0\) or every coroot is primitive.
\end{theorem}

\begin{proof}
This is proved in Proposition~\ref{prop:chevalley-criterion}.
\end{proof}

\subsection{Earlier methods}\label{subsec:previous}

Veldkamp proves the theorem for split semisimple groups over
perfect fields of characteristic greater than the largest Coxeter
number \cite[Prop.~2.1 and Thm.~3.1]{Veldkamp}.

Over an algebraically closed field \(k\) of characteristic \(p>0\),
Kac--Weisfeiler prove the Harish-Chandra isomorphism for almost simple
groups, excluding characteristic \(2\) for groups with nonprimitive
coroots \cite[Thm.~1 bis]{KW}. Their proof uses the \(p\)-centers
\[
Z_p(\mathfrak a)=k[x^p-x^{[p]}:x\in\mathfrak a],
\qquad \mathfrak a\in\{\fg_k,\fh_k\},
\]
where \(x^{[p]}\) is the restricted Lie-algebra power.
Via the descriptions of the \(p\)-centers as coordinate algebras of
Frobenius twists, Chevalley restriction identifies the subalgebras
\(Z_p(\fg_k)^{G_k}\subset U_k^{G_k}\) and
\(Z_p(\fh_k)^{W,\bullet}\subset S_k^{W,\bullet}\)
under Harish-Chandra projection. They then prove the full
Harish-Chandra isomorphism by comparing its source and target as
algebras over these identified invariant subalgebras of the two \(p\)-centers
\cite[Thm.~4(i), \S\S4.1--4.6, 5.3--5.4 and 8.1]{KW}.
Jantzen observes that their argument extends to reductive groups
with simply connected derived subgroup, without a good-prime
assumption \cite[\S9.3, Remark]{Jantzen}.
Polo--Tilouine further note that the method extends to general
reductive groups in positive characteristic, excluding characteristic
\(2\) for groups with nonprimitive coroots
\cite[\S2.6.2, Remark]{PT}.

\subsection{Strategy of proof}\label{subsec:our-method}

To prove the Harish-Chandra isomorphism over a general ring we cannot argue directly with \(p\)-centers. The main observation of this paper is the following argument, which easily reduces the Harish-Chandra isomorphism for \(G\) to the isomorphism for the Levi subgroups \(L_i\) of its minimal standard parabolic subgroups \(Q_i\).

Suppose that every coroot is primitive. For a standard parabolic \(Q\), put
\[
V_Q(R)=U_R/U_R\fu_{Q,R},\qquad
\cA_Q(R)=\operatorname{Ind}_{Q_R}^{G_R}V_Q(R),
\]
where \(\fu_Q=\Lie(R_u(Q))\) and induction means
\(\operatorname{Ind}_{Q_R}^{G_R}M=(\cO(G_R)\otimes_R M)^{Q_R}\) with the adjoint \(Q_R\)-action on \(V_Q(R)\).

The quotient maps \(V_Q(R)\to V_B(R)\) give strict filtered
\(G_R\)-equivariant embeddings
\(r_{Q,B}:\cA_Q(R)\into\cA_B(R)\).
With these identifications, Theorem~\ref{thm:parabolic-reconstruction}, which we call \emph{the parabolic intersection lemma}, states
\begin{equation}\label{eq:parabolic-reconstruction-intro}
U_R=\bigcap_i\cA_{Q_i}(R)\subset\cA_B(R)
\end{equation}
as filtered \(G_R\)-submodules.

It is not difficult to see that
\(\cA_{Q_i}(R)^{G_R}\simeq U_R(\fl_{i,R})^{L_{i,R}}\) and
\(\cA_B(R)^{G_R}\simeq S_R\), where \(\fl_i=\Lie(L_i)\).
Taking \(G_R\)-invariants in (\ref{eq:parabolic-reconstruction-intro}), and using the rank-one isomorphisms, gives
\begin{equation}\label{eq:hc_introduction}
\HC_R(U_R^{G_R})=\bigcap_i S_R^{s_i,\bullet}=S_R^{W,\bullet},
\end{equation}
where \(s_i\) is the simple reflection corresponding to \(Q_i\).

\subsection{Structure of the paper}
Section~\ref{sec:preliminaries} fixes conventions and notation.
Section~\ref{sec:reconstruction} proves the parabolic intersection lemma over an arbitrary ring. Section~\ref{sec:rankone} performs the rank-one calculation, proving
the Harish-Chandra isomorphism for groups with primitive coroots.
Section~\ref{sec:completion} constructs a central extension which resolves the nonprimitive coroots
and completes the proof. Its final
subsection establishes the criterion for Chevalley restriction over
arbitrary rings.

\subsection{Disclosure of AI use}

The author posed the research question and formulated the conjecture addressed in this paper. OpenAI's GPT-6 Astra generated an initial version of the proofs of the main theorems. The author then revised and rewrote the manuscript, including the proofs. During this work, the author consulted OpenAI's GPT-6.1 Sol. These tools also assisted with language editing, proofreading, and routine \LaTeX{} formatting. The author takes full responsibility for the mathematical arguments, references, and final content of the paper.

\section{Preliminaries}\label{sec:preliminaries}

\subsection{Root data and root vectors}

We retain the notation of the introduction. All group schemes and lattices without a coefficient subscript are defined over \(\Z\). We also regard a finite free \(R\)-module \(V\) as its associated vector bundle; with \(V^*=\operatorname{Hom}_R(V,R)\), this gives \(\cO(V^*)=\Sym_R(V)\). Reductivity means smoothness and affineness with connected reductive geometric fibers. Borel subgroups are understood fiberwise, and \(R_u(B)\) denotes the smooth normal unipotent radical.

Let
\[
X=X^*(T),\qquad Y=X_*(T),\qquad
\Phi\subset X,\qquad
\fB=\Lie(B),\qquad \fn=\Lie(R_u(B)).
\]
Differentiation identifies \(Y\) with \(\fh\), and hence \(Y\otimes_\Z R\) with \(\fh_R\). Write \(h_\alpha=\alpha^\vee\in Y\). Choose \(\Phi^+\) so that \(\fn\) is the sum of the positive root spaces. If \(B^-\) is the opposite Borel relative to \(T\) and \(\fn^-=\Lie(R_u(B^-))\), then
\[
\fg=\fn^-\oplus\fh\oplus\fn,\qquad
\fn^\pm=\bigoplus_{\alpha\in\Phi^+}\fg_{\pm\alpha}.
\]
Here \(\fg_\alpha\) is the root submodule of weight \(\alpha\), a free \(\Z\)-module of rank one. Choose opposite root vectors dual under the canonical perfect pairing of root modules, and let \(x_\alpha:\mathbb G_a\iso U_\alpha\) be the corresponding parametrizations, so that \(e_\alpha=d x_\alpha(1)\). Then \cite[Thm.~4.2.6 and Cor.~5.1.12]{Conrad} gives
\[
[e_\alpha,e_{-\alpha}]=h_\alpha.
\]
The same corollary gives the rank-one identities
\[
\begin{aligned}
\Ad(x_\alpha(t))e_{-\alpha}
&=e_{-\alpha}+t h_\alpha-t^2e_\alpha,\\
\Ad(x_\alpha(t))H&=H-t\alpha(H)e_\alpha
\qquad(H\in Y).
\end{aligned}
\]
These integral identities remain valid after base change.

\subsection{Actions, invariants, and weights}\label{subsec:actions}

For an affine group scheme \(H\) over \(R\), a rational \(H\)-module
is an \(R\)-module \(M\) with a coaction
\[
\delta_M:M\longrightarrow M\otimes_R\cO(H)
\]
satisfying the coassociativity and counit identities. We write
\[
M^H=\ker\bigl(\delta_M-(m\mapsto m\otimes1)\bigr).
\]
Thus invariance means fixedness under \(H(A)\) after every scalar
extension \(R\to A\). All action formulas below are understood in this
sense. For \(W\), invariants are simply the elements fixed by every
\(w\in W\).

The adjoint action of \(G_R\) on \(\fg_R\) extends to \(U_R\) and
\(\Sym_R(\fg_R)\), preserving their filtrations.

For the split torus, \(\cO(T_R)=R[X]\), with basis \(t^\lambda\)
indexed by the characters \(\lambda\in X\). Every rational
\(T_R\)-module has the weight decomposition
\begin{equation}\label{eq:weights}
M=\bigoplus_{\lambda\in X}M_\lambda,\qquad
M_\lambda=\{m:\delta_M(m)=m\otimes t^\lambda\},\qquad
M^{T_R}=M_0.
\end{equation}
This follows by comparing coefficients of the basis \(t^\lambda\) in
the coaction identities. Weights lie in the character lattice \(X\),
and the decomposition requires no flatness assumption on \(M\).

\subsection{The dot action}

The dot action of \(W\) on \(\fh^*\) is
\[
w\bullet\lambda=w(\lambda+\rho)-\rho,
\qquad
\rho=\frac12\sum_{\alpha>0}\alpha.
\]
It is defined over \(\Z\), since
\[
w\rho-\rho
=-\sum_{\substack{\alpha\in\Phi^+\\w^{-1}\alpha\in\Phi^-}}\alpha
\]
belongs to the root lattice. We therefore do not need to assume that \(\rho\in X\) to define the dot action.

We identify the symmetric algebra \(S_R=\Sym_R(\fh_R)\) with \(\cO(\fh_R^*)\), where \(H\in Y\) represents the linear function \(\lambda\mapsto\lambda(H)\). The dot action on functions \(f\in\cO(\fh_R^*)\) is the pullback action
\[
(w\bullet f)(\lambda)=f(w^{-1}\bullet\lambda).
\]
On the linear generators of the symmetric algebra \(S_R\), a simple reflection acts by
\begin{equation}\label{eq:dot}
s_\alpha\bullet H=H-\alpha(H)(h_\alpha+1)
\qquad(H\in Y).
\end{equation}
Here \(1\) denotes the constant polynomial. Thus \(s_\alpha\bullet H=s_\alpha H-\alpha(H)\), where \(s_\alpha H=H-\alpha(H)h_\alpha\) is the ordinary Weyl group action. Therefore, the dot action preserves the polynomial-degree filtration
\[
F_dS_R=\bigoplus_{j=0}^d\Sym_R^j(\fh_R).
\]
Moreover, for \(f\in F_dS_R\),
\[
w\bullet f\equiv wf\pmod{F_{d-1}S_R},
\]
where \(wf\) denotes the ordinary Weyl action and \(F_{-1}S_R=0\). Thus the induced action on \(\gr S_R\simeq S_R\) is the ordinary Weyl group action.

\subsection{An exactness lemma over \texorpdfstring{\(\Z\)}{Z}}

Note that a \(\Z\)-module is flat if and only if it is torsion-free, with no finite-generation hypothesis \cite[Tag 0AUW(1)]{Stacks}. We will use the following criterion throughout the paper.

\begin{lemma}\label{lem:flat}\leavevmode
\begin{enumerate}
\item\label{item:flat-injection}
Let \(j:A\to M\) be a homomorphism of flat \(\Z\)-modules. If \(j\otimes k\) is injective for \(k=\Q\) and every \(k=\F_p\), then \(j\) is injective with flat cokernel. In particular, \(j\otimes L\) is injective for every \(\Z\)-module \(L\).
\item\label{item:flat-complex}
Let \(A\xrightarrow{j}M\xrightarrow{D}N\) be a complex of flat \(\Z\)-modules. Suppose
\[
0\longrightarrow A\otimes k
\xrightarrow{j_k}M\otimes k
\xrightarrow{D_k}N\otimes k
\]
is exact at its first two terms for \(k=\Q\) and every \(k=\F_p\). Then it is exact at those terms after tensoring with any \(\Z\)-module. Moreover, both \(M/j(A)\) and \(\operatorname{coker}D\) are flat over \(\Z\).
\end{enumerate}
\end{lemma}
\begin{proof}
For~\eqref{item:flat-injection}, torsion-freeness of \(A\) and injectivity of \(j\otimes\Q\) imply that \(j\) is injective. Set \(Q=M/j(A)\), and write \([m]=m+j(A)\) for the class of \(m\in M\) in \(Q\). Fix a prime \(p\), and suppose that \(p[m]=0\). Then \(pm=j(a)\) for some \(a\in A\). Reducing modulo \(p\) and using injectivity of \(j_{\F_p}\), we obtain \(a=pa'\) for some \(a'\in A\). Hence
\[
p(m-j(a'))=0.
\]
Since \(M\) is torsion-free, \(m=j(a')\), so \([m]=0\). Thus \(Q\) has no \(p\)-torsion for any prime \(p\), and is therefore torsion-free, hence flat. It follows that the exact sequence \(0\to A\to M\to Q\to0\) remains exact after tensoring with any \(\Z\)-module \(L\), proving the last assertion.

For~\eqref{item:flat-complex}, part~\eqref{item:flat-injection} shows that \(j\) is injective and \(Q=M/j(A)\) is flat. Since \(D\circ j=0\), the map \(D\) factors through
\[
\overline D:Q\longrightarrow N,\qquad [m]\longmapsto D(m).
\]
Exactness over \(\Q\) implies that \(\overline D\otimes\Q\) is injective; since \(Q\) is torsion-free, \(\overline D\) is injective as well. As \(Q\) is the cokernel of \(j\), reduction modulo \(p\) gives the canonical identifications
\[
Q/pQ\simeq M/(j(A)+pM)
\simeq (M/pM)/j_{\F_p}(A/pA).
\]
Under these identifications, exactness at \(M/pM\) says precisely that \(\overline D\otimes\F_p\) is injective for every prime \(p\). Applying part~\eqref{item:flat-injection} to \(\overline D:Q\to N\) now shows that \(C=\operatorname{coker}\overline D=\operatorname{coker}D\) is flat. Since \(Q\) and \(C\) are flat, both short exact sequences
\[
0\to A\to M\to Q\to0,\qquad
0\to Q\to N\to C\to0
\]
remain exact after tensoring with any \(\Z\)-module \(L\). Combining them proves the required exactness at \(A\otimes L\) and \(M\otimes L\).
\end{proof}

\subsection{PBW theorem and Harish-Chandra projection}

PBW gives
\[
U_R\simeq U_R(\fn_R^-)\otimes_R S_R\otimes_R U_R(\fn_R)
\]
as \(R\)-modules, with negative root vectors ordered first, then a basis of \(Y\), then positive root vectors. It also identifies \(\gr U_R\) with \(\Sym_R(\fg_R)\); see \cite[Ch.~I, \S2, no.~7, Thm.~1]{Bourbaki123}, or the explicit arbitrary-ring statement in \cite[I.1.40]{Milne}. Let \(\varepsilon_\pm:U_R(\fn_R^\pm)\to R\) be the augmentation maps, which send \(1\) to \(1\) and every element of \(\fn_R^\pm\) to zero. Under the PBW identification, the \(R\)-linear Harish-Chandra projection is
\[
\HC_R=\varepsilon_-\otimes\operatorname{id}_{S_R}\otimes\varepsilon_+:
U_R\longrightarrow S_R.
\]

The following multiplicativity statement over arbitrary rings is also recorded in \cite[\S2.6.1]{PT}.

\begin{lemma}\label{lem:multiplicativity}
Harish-Chandra projection is multiplicative on \(U_R^{T_R}\), and hence on \(U_R^{G_R}\).
\end{lemma}
\begin{proof}
Let \(V_B=U_R/U_R\fn_R\), the quotient by the left ideal generated by \(\fn_R\). It has the left \(U_R\)-action by multiplication and the \(B_R\)-action induced by the adjoint action. Write \(\bar1\) for the class of \(1\). PBW gives
\[
V_B\simeq U_R(\fn_R^-)\otimes_R S_R
\]
and an injection \(\iota:S_R\into V_B\), \(s\mapsto s\cdot\bar1\). Every PBW monomial involving \(\fn_R^-\) has weight in the nonzero negative root cone. By~\eqref{eq:weights}, \(V_B^{T_R}=\iota(S_R)\). Thus the \(T_R\)-equivariant quotient map \(q:U_R\to V_B\) satisfies
\[
q(a)=\iota(\HC_R(a))=\HC_R(a)\cdot\bar1
\qquad(a\in U_R^{T_R}).
\]
Differentiating the torus action shows that \(a\) commutes with \(\fh_R\), and hence with \(S_R\). For \(a,b\in U_R^{T_R}\), the left \(U_R\)-module structure on \(V_B\) therefore gives
\[
\begin{aligned}
q(ab)&=a\cdot q(b)
=a\cdot\bigl(\HC_R(b)\cdot\bar1\bigr)\\
&=\HC_R(b)\cdot\bigl(a\cdot\bar1\bigr)
=\iota\bigl(\HC_R(b)\HC_R(a)\bigr).
\end{aligned}
\]
Since \(ab\) is also \(T_R\)-fixed and \(\iota\) is injective, this proves multiplicativity.
\end{proof}

\subsection{Torsion-free base ring}

\begin{lemma}\label{lem:center}
For every \(R\), one has \(U_R^{G_R}\subseteq Z(U_R)\). Equality holds if and only if \(R\) is torsion-free as a \(\Z\)-module.
\end{lemma}
\begin{proof}
Differentiating the adjoint action shows that a group-invariant element commutes with \(\fg_R\). Since \(\fg_R\) generates \(U_R\), the element is central.

Assume that \(R\) is torsion-free, and let \(z\in Z(U_R)\). By PBW, \(U_R\) is free over \(R\) and hence torsion-free over \(\Z\). Write the torus weight decomposition as \(z=\sum_\lambda z_\lambda\). For \(H\in Y=\fh\),
\[
0=[H,z]=\sum_\lambda\lambda(H)z_\lambda.
\]
If \(\lambda\ne0\), some \(H\in Y\) has \(\lambda(H)\ne0\). Torsion-freeness therefore gives \(z_\lambda=0\) for every \(\lambda\ne0\), so \(z=z_0\) is \(T_R\)-fixed.

For a root subgroup, let \(\rho_\alpha(t)=\Ad(x_\alpha(t))\) on \(U_R\), and write
\[
\rho_\alpha(t)z=\sum_{m\ge0}v_mt^m\in U_R[t].
\]
Differentiate the polynomial identity
\(\rho_\alpha(t+s)z=\rho_\alpha(t)\rho_\alpha(s)z\) with respect to \(s\) at \(s=0\). Since the differential of the root-subgroup action is \(\operatorname{ad}(e_\alpha)\), this gives
\[
\frac{d}{dt}\rho_\alpha(t)z
=\rho_\alpha(t)[e_\alpha,z]=0.
\]
Thus \(m v_m=0\) for \(m>0\); torsion-freeness implies \(v_m=0\), and \(z\) is fixed by every root subgroup.

By \cite[Thm.~5.1.13]{Conrad}, multiplication identifies each unipotent radical with the ordered product of its root subgroups, and
\[
\Omega:=R_u(B^-)\times T\times R_u(B)\longrightarrow G
\]
is an open immersion with fiberwise dense image, the \emph{big cell}. Thus the coaction identity for \(z\) holds on \(\Omega_R\). For \(k=\Q\) and every \(k=\F_p\), the group \(G_k\) is smooth, hence reduced, and \(\Omega_k\) is dense, so restriction \(\cO(G_k)\to\cO(\Omega_k)\) is injective. Both \(\cO(G)\) and \(\cO(\Omega)\) are \(\Z\)-flat, so Lemma~\ref{lem:flat}\eqref{item:flat-injection} gives injectivity of \(\cO(G_R)\to\cO(\Omega_R)\) for every \(R\). Tensoring with the free \(R\)-module \(U_R\) now shows that the coaction identity holds on \(G_R\), proving \(z\in U_R^{G_R}\).

Conversely, suppose that \(R\) has nonzero additive torsion. Choose a prime \(p\) and \(0\ne a\in R\) with \(pa=0\), and choose a root \(\alpha\). We first establish the integral divisibility
\begin{equation}\label{eq:root-p-divisibility}
[e_\alpha^p,\fg]\subseteq pU_{\Z}.
\end{equation}
Write the root-subgroup action on \(\fg\) as
\[
\Ad(x_\alpha(t))\big|_{\fg}
=\sum_{m\ge0}A_mt^m,\qquad
A_m\in\operatorname{End}_{\Z}(\fg).
\]
The additive group law gives
\[
A_iA_j=\binom{i+j}{i}A_{i+j},\qquad
A_1=\operatorname{ad}(e_\alpha).
\]
Consequently \(\operatorname{ad}(e_\alpha)^p=p!A_p\), so reduction modulo \(p\) gives
\(\operatorname{ad}(\overline e_\alpha)^p=0\) on \(\fg_{\F_p}\). In an associative algebra of characteristic \(p\), the binomial identity gives
\[
\operatorname{ad}(\overline e_\alpha^{\,p})
=\operatorname{ad}(\overline e_\alpha)^p.
\]
Thus \(\overline e_\alpha^{\,p}\) commutes with the generators of \(U_{\F_p}\). Since \(U_{\Z}/pU_{\Z}\simeq U_{\F_p}\), this proves~\eqref{eq:root-p-divisibility}.

Now \(z=a e_\alpha^p\in U_R\) commutes with \(\fg_R\) by~\eqref{eq:root-p-divisibility} and \(pa=0\), so it is central. PBW shows that \(z\ne0\). Its torus weight is \(p\alpha\ne0\) in the character lattice \(X\), so the weight decomposition~\eqref{eq:weights} shows that \(z\) is not \(T_R\)-fixed, and hence not \(G_R\)-fixed. This proves the converse.
\end{proof}

\section{The parabolic intersection lemma}\label{sec:reconstruction}

We prove the parabolic intersection lemma in this section.

\subsection{The modules \texorpdfstring{\(\mathcal C_Q(R)\) and \(\cA_Q(R)\)}{C\_Q(R) and A\_Q(R)}}\label{subsec:classical-comparisons}

Let \(Q=LU_Q\) be a standard parabolic, with Levi subgroup \(L\)
containing \(T\). Write \(\fl=\Lie(L)\) and
\(\fu_Q=\Lie(U_Q)\), and let \(\fu_Q^-\) be the opposite unipotent
Lie algebra. The left ideal \(U_R\fu_{Q,R}\) is \(Q_R\)-stable, so
\[
V_Q(R)=U_R/U_R\fu_{Q,R}
\]
inherits the adjoint \(Q_R\)-action and the quotient PBW filtration. We form the induced module
\[
\cA_Q(R)=\operatorname{Ind}_{Q_R}^{G_R}V_Q(R),
\qquad
\operatorname{Ind}_{Q_R}^{G_R}M
=\bigl(\cO(G_R)\otimes_R M\bigr)^{Q_R}.
\]
A tensor \(F\) is viewed as a regular \(M\)-valued function on \(G_R\).
The invariance condition and the left \(G_R\)-action are
\begin{equation}\label{eq:induced-actions}
F(gq)=q^{-1}\cdot F(g),\qquad
(a\cdot F)(g)=F(a^{-1}g).
\end{equation}
Here \(Q_R\) acts on the tensor product by right translation on
\(\cO(G_R)\) and its given action on \(M\). The filtration is
\[
F_d\cA_Q(R)=
\bigl(\cO(G_R)\otimes_R F_dV_Q(R)\bigr)^{Q_R}.
\]
It is exhaustive, with \(F_{-1}\cA_Q(R)=0\).

Write \(q_Q:U_R\to V_Q(R)\) for the quotient. Since
\(\fu_Q\subseteq\fn\), there is also a \(B_R\)-equivariant quotient
\(q_{Q,B}:V_Q(R)\to V_B(R)\). They give filtered
\(G_R\)-equivariant maps
\begin{equation}\label{eq:filtered-comparison-maps}
\begin{aligned}
j_{Q,R}:U_R&\longrightarrow\cA_Q(R),
&j_{Q,R}(u)(g)&=q_Q(g^{-1}\cdot u),\\
r_{Q,B}:\cA_Q(R)&\longrightarrow\cA_B(R),
&r_{Q,B}(F)(g)&=q_{Q,B}(F(g)).
\end{aligned}
\end{equation}
The equivariance condition in~\eqref{eq:induced-actions} verifies
these formulas, and \(r_{Q,B}j_{Q,R}=j_{B,R}\).
For \(Q=G\), evaluation at the identity identifies
\(\cA_G(R)\) with \(U_R\), with inverse \(j_{G,R}\).
We prove below that all these maps are strict embeddings.

The classical counterpart of \(\cA_Q(R)\) is
\[
\mathcal C_Q(R)=
\bigl(\cO(G_R)\otimes_R\Sym_R(\fg_R/\fu_{Q,R})\bigr)^{Q_R}.
\]
PBW identifies
\(\gr V_Q(R)=\Sym_R(\fg_R/\fu_{Q,R})\).
Taking leading PBW symbols gives a canonical graded injection
\[
\iota_Q:\gr\cA_Q(R)\into\mathcal C_Q(R).
\]
To see injectivity, let \(F\in F_d\cA_Q(R)\) have zero symbol in degree \(d\).
Since \(\cO(G_R)\) is flat over \(R\), this means that
\(F\in\cO(G_R)\otimes_R F_{d-1}V_Q(R)\).
The tensor \(F\) is \(Q_R\)-invariant, so it belongs to
\(F_{d-1}\cA_Q(R)\), as required.
The grading on \(\mathcal C_Q(R)\) is given by the polynomial degree
in the symmetric algebra.

The classical counterparts of \(j_{Q,R}\) and \(r_{Q,B}\) are obtained
from~\eqref{eq:filtered-comparison-maps} by replacing \(q_Q\) and
\(q_{Q,B}\) with their associated graded maps. We denote them by
\[
\begin{aligned}
j^{\rm cl}_{Q,R}&:\Sym_R(\fg_R)\longrightarrow\mathcal C_Q(R),\\
r^{\rm cl}_{Q,B}&:\mathcal C_Q(R)\longrightarrow\mathcal C_B(R).
\end{aligned}
\]
These maps are graded and \(G_R\)-equivariant. By construction,
\begin{equation}\label{eq:graded-maps}
\iota_Q\circ\gr(j_{Q,R})=j^{\rm cl}_{Q,R},\qquad
\iota_B\circ\gr(r_{Q,B})=r^{\rm cl}_{Q,B}\circ\iota_Q.
\end{equation}

To interpret \(\mathcal C_Q(R)\) geometrically, put
\(E_Q=(\fg/\fu_Q)^*\subset\fg^*\), the annihilator of \(\fu_Q\).
For any commutative ring \(A\), the coadjoint action of \(G(A)\) on
\(\fg_A^*=\operatorname{Hom}_A(\fg_A,A)\) is given by
\[
(g\xi)(x)=\xi(\Ad(g^{-1})x)
\qquad(g\in G(A),\ \xi\in\fg_A^*,\ x\in\fg_A).
\]
Restricting this action gives a \(Q\)-action on \(E_Q\).
Let \(G\times^Q E_Q\) be the quotient of \(G\times E_Q\) by the right
\(Q\)-action \((g,\xi)\cdot q=(gq,q^{-1}\xi)\).
We write \([g,\xi]\) for the class represented by \((g,\xi)\). 

Moreover, \(G\times^Q E_Q\) is a vector bundle over \(G/Q\), with
bundle projection \([g,\xi]\mapsto gQ\). Its total space identifies with
the incidence scheme
\[
X_Q=\{(\chi,gQ)\in\fg^*\times G/Q:
       \chi(\Ad(g)\fu_Q)=0\}.
\]
The identification is given by
\[
[g,\xi]\longmapsto(g\xi,gQ).
\]
For \(Q=G\), one has \(X_G=\fg^*\).
Descent of regular functions along the \(Q_R\)-torsor
\(G_R\times E_{Q,R}\to X_{Q,R}\) gives
\begin{equation}\label{eq:incidence-functions}
\mathcal C_Q(R)=\Gamma(X_{Q,R},\cO)
\end{equation}
for every \(R\) \cite[Tag 023Q]{Stacks}.

The incidence map
\[
\mu_Q:X_Q\longrightarrow\fg^*,\qquad
(\chi,gQ)\longmapsto\chi
\]
is proper: the incidence condition is closed, and \(G/Q\) is proper
\cite[Cor.~5.2.8]{Conrad}.
The scheme \(X_Q\) is smooth over \(\Z\), and its geometric fibers
are integral of dimension
\[
\dim(G/Q)+\dim(\fg/\fu_Q)=\dim G.
\]
The action of \(G\) on \(X_Q\) is given by \(a(\chi,gQ)=(a\chi,agQ)\).
For \(B\subset Q\), forgetting part of the flag gives
\begin{equation}\label{eq:incidence-projections}
\pi_Q:X_B\longrightarrow X_Q,\qquad
(\chi,gB)\longmapsto(\chi,gQ),\qquad
\mu_Q\pi_Q=\mu_B.
\end{equation}
The geometric description identifies \(j^{\rm cl}_{Q,R}\) and \(r^{\rm cl}_{Q,B}\) as pullbacks.
Since \(\gr(q_Q)\) restricts polynomials on \(\fg_R^*\) to \(E_{Q,R}\),
we have
\[
j^{\rm cl}_{Q,R}(f)([g,\xi])=(g^{-1}\cdot f)(\xi)=f(g\xi)
\qquad(f\in\Sym_R(\fg_R)).
\]
Thus \(j^{\rm cl}_{Q,R}=\mu_{Q,R}^*\).
Similarly, \(\gr(q_{Q,B})\) restricts polynomials from \(E_{Q,R}\)
to \(E_{B,R}\). The map \(\pi_{Q,R}\) sends \([g,\xi]\in X_{B,R}\)
to \([g,\xi]\in X_{Q,R}\), so \(r^{\rm cl}_{Q,B}=\pi_{Q,R}^*\).

\begin{lemma}\label{lem:coadjoint-density}
For every algebraically closed field \(k\), the proper incidence map
\[
\mu_B:X_{B,k}\longrightarrow\fg_k^*
\]
is surjective and generically finite.
For every standard parabolic \(Q\), the map
\(\pi_Q:X_{B,k}\to X_{Q,k}\) is surjective as well.
\end{lemma}
\begin{proof}
We use the Bruhat decomposition argument of \cite[Lem.~2.5]{PT}.
Choose \(\chi_0\in\fg_k^*\)
vanishing on \(\fB_k\) and the negative nonsimple root vectors, with
\[
\chi_0(e_{-\alpha})=1
\quad(\alpha\text{ simple}).
\]
A point \(gB\) lies in the incidence fiber over \(\chi_0\) exactly
when \(\chi_0(\Ad(g)\fn_k)=0\).
Write \(g=u n_w b\) by Bruhat decomposition
\cite[Thm.~1.4.12]{Conrad}, with
\(u\in R_u(B)(k)\), \(b\in B(k)\), and \(n_w\) representing \(w\in W\).
If \(w\ne1\), choose a simple \(\alpha\) with \(w^{-1}\alpha<0\).
Then \(e_{-\alpha}\in\Ad(n_w)\fn_k\), so the incidence condition gives
\[
0=\chi_0(\Ad(u)e_{-\alpha})=\chi_0(e_{-\alpha})=1.
\]
For the middle equality, positive root groups change
\(e_{-\alpha}\) only by elements of \(\fB_k\): any weight
\(-\alpha+j\gamma\) that occurs, with \(\gamma>0\) and \(j>0\),
is zero or a positive root.

Thus this fiber has the single geometric point \(B\), and hence
dimension zero. The image of \(\mu_B\) is closed by properness.
The fiber-dimension inequality and
\(\dim X_{B,k}=\dim\fg_k^*\) imply that this closed image has full
dimension, so it is all of \(\fg_k^*\).
The equality of dimensions also shows that \(\mu_B\) is generically finite.

To prove surjectivity of \(\pi_Q\), let \((\chi,gQ)\in X_{Q,k}\)
and put \(\xi=g^{-1}\chi\in E_{Q,k}\).
Apply the surjectivity just proved to the Levi \(L\), with Borel
\(B\cap L\), and to the restriction \(\xi|_{\fl_k}\).
It gives \(\ell\in L(k)\) such that \(\xi\) annihilates
\(\Ad(\ell)(\fn_k\cap\fl_k)\).
It already annihilates the \(L\)-stable space \(\fu_{Q,k}\).
Since
\(\fn_k=\fu_{Q,k}\oplus(\fn_k\cap\fl_k)\),
the point \((\chi,g\ell B)\) belongs to \(X_{B,k}\) and lifts
\((\chi,gQ)\).
\end{proof}

\begin{lemma}\label{lem:comparisons}
For every split connected reductive \(G\), standard parabolic \(Q\),
and ring \(R\), the maps \(r_{Q,B}\) and \(j_{Q,R}\) are strict
filtered \(G_R\)-equivariant embeddings.
Their classical analogues \(r^{\rm cl}_{Q,B}=\pi_{Q,R}^*\) and
\(j^{\rm cl}_{Q,R}=\mu_{Q,R}^*\) are graded
\(G_R\)-equivariant embeddings.
\end{lemma}
\begin{proof}
Lemma~\ref{lem:coadjoint-density} and
\(\mu_Q\pi_Q=\mu_B\) show that both \(\pi_Q\) and \(\mu_Q\)
are surjective on every geometric fiber over \(\Z\).

We check that their pullbacks remain injective after arbitrary base
change. Let \(f\) denote either map, choose an affine open \(V\)
in its target, and cover \(f^{-1}(V)\) by finitely many affine opens
\(U_j\). All these schemes are smooth over \(\Z\), so their
coordinate rings are \(\Z\)-flat.
For \(k=\Q\) and every \(k=\F_p\), the map
\[
\cO(V)\otimes k\longrightarrow
\bigoplus_j\cO(U_j)\otimes k
\]
is injective: after extending to an algebraic closure, \(f_k\)
is surjective and \(V_k\) is reduced.
Lemma~\ref{lem:flat}\eqref{item:flat-injection} makes this map
injective after tensoring with every \(R\).
Working locally on the target proves that \(f_R^*\) is injective
on global functions. This proves injectivity of \(j^{\rm cl}_{Q,R}\) and \(r^{\rm cl}_{Q,B}\).

Equation~\eqref{eq:graded-maps} and the injections \(\iota_Q\) now show
that \(\gr(j_{Q,R})\) and \(\gr(r_{Q,B})\) are injective.
An injective associated graded map preserves the degree of every
nonzero element, so the corresponding filtered maps are injective and strict.
Equivariance and the grading follow from the constructions.
\end{proof}

From now on, we identify \(\cA_Q(R)\) and \(\mathcal C_Q(R)\)
with their images in \(\cA_B(R)\) and \(\mathcal C_B(R)\)
under \(r_{Q,B}\) and \(r^{\rm cl}_{Q,B}\), respectively.
These identifications respect the filtrations, gradings, and \(G_R\)-actions.
For \(Q=G\), this also identifies \(U_R\) and \(\Sym_R(\fg_R)\)
with their images under \(j_{B,R}\) and \(j^{\rm cl}_{B,R}\), respectively.

\subsection{Proof of the lemma}

Let \(Q_1,\ldots,Q_r\) be the minimal standard parabolics properly
containing \(B\), with simple reflections \(s_i\), and write
\(\pi_i=\pi_{Q_i}\).
We first prove the parabolic intersection lemma for \(\mathcal C_B(R)\).

Over an algebraically closed field \(k\), put
\[
P=\Sym_k(\fg_k),\qquad K=\Frac(P),\qquad K'=k(X_{B,k}).
\]
The dominant, generically finite map \(\mu_B\) identifies \(K\)
with a subfield of \(K'\), and \(\mathcal C_B(k)\) with a subalgebra
of \(K'\).

\begin{lemma}\label{lem:polynomial-descent}
Inside \(K'\), one has
\[
\mathcal C_B(k)\cap K=P.
\]
\end{lemma}
\begin{proof}
Let \(F=a/b\in\mathcal C_B(k)\cap K\), with coprime \(a,b\in P\).
If an irreducible \(q\) divides \(b\), surjectivity of \(\mu_B\)
gives a point \(x\in X_{B,k}\) over the generic point of the divisor
\((q)\). In the local ring \(\cO_{X_{B,k},x}\), the element \(a\)
is a unit and \(b\) is a nonunit. This contradicts \(a=bF\), since
\(F\) is regular. Thus \(b\) is a unit and \(F\in P\).
The reverse inclusion is pullback along \(\mu_B\).
\end{proof}

\begin{lemma}\label{lem:bruhat-compatibility}
Let \(k\) be algebraically closed.
Any two \(k\)-points of a fiber of \(\mu_B:X_{B,k}\to\fg_k^*\) can be joined
by a sequence of \(k\)-points in that fiber such that consecutive points
have the same image under some \(\pi_i\).
\end{lemma}
\begin{proof}
Let the points be \((\chi,gB)\) and \((\chi,g'B)\).
Write \(g^{-1}g'=b n_w b'\) by Bruhat decomposition and put \(a=gb\).
Then \(gB=aB\) and \(g'B=a n_wB\).
Thus \(\eta=a^{-1}\chi\) annihilates both
\(\fn_k\) and \(\Ad(n_w)\fn_k\).

Choose a reduced expression \(w=s_{i_1}\cdots s_{i_m}\), put
\(v_0=1\), \(v_j=s_{i_1}\cdots s_{i_j}\), and choose representatives
as products of the simple-reflection representatives.
For each prefix,
\[
v_j\Phi^+\subseteq\Phi^+\cup w\Phi^+.
\]
Indeed, \(v_j\Phi^+\) is obtained from \(v_{j-1}\Phi^+\) by replacing
\(v_{j-1}\alpha_{i_j}\) with its negative. Reducedness ensures that
every removed root is positive, so any negative root in \(v_j\Phi^+\)
remains in \(w\Phi^+\).
Consequently
\[
\Ad(n_{v_j})\fn_k\subseteq\fn_k+\Ad(n_w)\fn_k.
\]
Hence \(\eta\) annihilates every \(\Ad(n_{v_j})\fn_k\), and the
points \((\chi,a n_{v_j}B)\) all lie in the original fiber.
Consecutive flags have the same image in \(G/Q_{i_j}\), since
\(n_{s_{i_j}}\in Q_{i_j}(k)\). Their incidence points therefore have
the same image under \(\pi_{i_j}\), as required.
\end{proof}

\begin{proposition}\label{prop:classical-reconstruction}
Assume that every coroot is primitive in \(Y\). For every commutative ring \(R\),
\[
\Sym_R(\fg_R)=\bigcap_i\mathcal C_{Q_i}(R)
\subset\mathcal C_B(R).
\]
\end{proposition}
\begin{proof}
First let \(k\) be algebraically closed.
Choose \(\lambda\in\fh_k^*\) with
\(\lambda(h_\alpha)\ne0\) for every root, and extend it by zero on
the root spaces. Such a choice exists because 
each \(h_\alpha\) is nonzero over \(k\).
Near the point \([1,\lambda]\), we trivialize \(X_{B,k}\) over the
big cell \(U_k^-B_k/B_k\subset G_k/B_k\), where \(U^-=R_u(B^-)\),
via \((u,\xi)\mapsto[u,\xi]\) for \((u,\xi)\in U_k^-\times E_{B,k}\).
In these coordinates, the differential of \(\mu_B\) at \([1,\lambda]\) is
\[
\fn_k^-\oplus E_{B,k}\longrightarrow\fg_k^*,\qquad
(v,\eta)\longmapsto\operatorname{ad}^*(v)\lambda+\eta.
\]
Modulo \(E_{B,k}\), the negative root vector \(e_{-\alpha}\)
maps to \(\lambda(h_\alpha)e_\alpha^*\), where \(e_\alpha^*\) is
the dual coordinate. Thus the differential is an isomorphism.
Since source and target are smooth of the same dimension,
\(\mu_B\) is generically \'{e}tale, and the finite extension
\(K'/K\) is separable \cite[Tags 02GH and 02NV]{Stacks}.

Let \(F\in\bigcap_i\mathcal C_{Q_i}(k)\). By
\eqref{eq:incidence-functions}, we view \(F\) as a regular function
on \(X_{B,k}\). For each \(i\), there is a regular function \(F_i\)
on \(X_{Q_i,k}\) such that \(F=\pi_i^*F_i\).
Let \(\overline K\) be an algebraic closure of \(K\).
After base change along \(k\hookrightarrow\overline K\), write
\(F_{\overline K}\) and \(F_{i,\overline K}\) for the resulting
regular functions on \(X_{B,\overline K}\) and
\(X_{Q_i,\overline K}\), respectively. They satisfy
\[
F_{\overline K}=\pi_{i,\overline K}^*F_{i,\overline K},
\qquad
\pi_{i,\overline K}:X_{B,\overline K}\longrightarrow X_{Q_i,\overline K}.
\]

The inclusion \(P\hookrightarrow K\hookrightarrow\overline K\)
determines an \(\overline K\)-point \(\chi_{\rm gen}\) of
\(\fg_{\overline K}^*\), whose image in \(\fg_k^*\) is the generic point.
The fiber of \(\mu_{B,\overline K}\) over \(\chi_{\rm gen}\) is
the geometric generic fiber of the original map \(\mu_{B,k}\),
namely \(\Spec(K'\otimes_K\overline K)\).

Apply Lemma~\ref{lem:bruhat-compatibility} with \(\overline K\) in place
of \(k\) to \(\mu_{B,\overline K}:X_{B,\overline K}\to\fg_{\overline K}^*\).
It says that any two \(\overline K\)-points of this fiber can be
joined by a chain of points in the fiber such that each consecutive
pair \(x,x'\) has the same image under some \(\pi_{i,\overline K}\).
For such a pair, the pullback identity gives
\[
F_{\overline K}(x)
=F_{i,\overline K}\bigl(\pi_{i,\overline K}(x)\bigr)
=F_{i,\overline K}\bigl(\pi_{i,\overline K}(x')\bigr)
=F_{\overline K}(x').
\]
Following the chain shows that \(F_{\overline K}\) takes the same
value at every \(\overline K\)-point of this fiber.

The \(\overline K\)-points of the geometric generic fiber
\(\Spec(K'\otimes_K\overline K)\) correspond to the \(K\)-embeddings
\(\sigma:K'\hookrightarrow\overline K\), via
\(a\otimes b\mapsto\sigma(a)b\).
The restriction of \(F_{\overline K}\) to the fiber is \(F\otimes1\),
so its value at the point corresponding to \(\sigma\) is \(\sigma(F)\).
Since the finite extension \(K'/K\) is separable and all
\(\sigma(F)\) are equal, we obtain \(F\in K\).
Lemma~\ref{lem:polynomial-descent} then gives \(F\in P\).
The converse inclusion \(P\subseteq\bigcap_i\mathcal C_{Q_i}(k)\)
is obvious from the identities \(\mu_B=\mu_{Q_i}\circ\pi_i\).

Since field extensions are faithfully flat, we deduce the assertion over every field from
the algebraically closed case. All maps appearing here preserve degree, so the equality
holds degree by degree.

To prove the equality over an arbitrary commutative ring \(R\),
we apply Lemma~\ref{lem:flat}\eqref{item:flat-complex} to a
complex of flat \(\Z\)-modules.
Fix \(d\ge0\) and put
\[
\begin{aligned}
M_{Q,d}&=\cO(G)\otimes_{\Z}\Sym_{\Z}^d(\fg/\fu_Q),\\
M_d&=\bigoplus_i M_{Q_i,d},\qquad
N_d=\bigoplus_i\bigl(M_{Q_i,d}\otimes\cO(Q_i)\bigr)
\oplus M_{B,d}^{\oplus(r-1)}.
\end{aligned}
\]
Let \(\delta_i\) be the diagonal \(Q_i\)-coaction, and let
\(b_i:M_{Q_i,d}\to M_{B,d}\) be restriction from \(E_{Q_i}\) to \(E_B\).
Define
\[
D_d((F_i))=\bigl((\delta_i(F_i)-F_i\otimes1)_i,
(b_i(F_i)-b_1(F_1))_{i>1}\bigr).
\]
The pullbacks along the \(\mu_{Q_i}\) give a complex of flat \(\Z\)-modules
\begin{equation}\label{eq:classical-reconstruction-complex}
0\longrightarrow\Sym_{\Z}^d(\fg)
\xrightarrow{j_d}M_d\xrightarrow{D_d}N_d.
\end{equation}

For \(k=\Q\) and every \(k=\F_p\),
Lemma~\ref{lem:comparisons} gives injectivity of \(j_d\otimes k\),
and the field case gives exactness of~\eqref{eq:classical-reconstruction-complex}
at \(M_d\otimes k\).
Lemma~\ref{lem:flat}\eqref{item:flat-complex} therefore gives exactness
after tensoring with every \(R\).
The kernel of \(D_d\otimes R\) consists of tuples of degree-\(d\)
elements of the \(\mathcal C_{Q_i}(R)\) that agree in \(\mathcal C_B(R)\).
Their common value gives the degree-\(d\) intersection, proving the
equality over \(R\).
\end{proof}

\begin{theorem}[The parabolic intersection lemma]\label{thm:parabolic-reconstruction}
Assume that every coroot is primitive in \(Y\). For every commutative ring \(R\),
\[
U_R=\bigcap_i\cA_{Q_i}(R)\subset\cA_B(R).
\]
\end{theorem}
\begin{proof}
The identities \(r_{Q_i,B}j_{Q_i,R}=j_{B,R}\) give
\(U_R\subseteq\bigcap_i\cA_{Q_i}(R)\).
For the converse, let \(F\in F_d\cA_B(R)\) lie in the intersection.
Then \(F\in F_d\cA_{Q_i}(R)\) for every \(i\), so~\eqref{eq:graded-maps} gives
\(\iota_B(\sigma_d(F))\in\bigcap_i\mathcal C_{Q_i}(R)\).
Proposition~\ref{prop:classical-reconstruction} shows that
\(f:=\iota_B(\sigma_d(F))\) belongs to \(\Sym_R^d(\fg_R)\).
Choose \(u\in F_dU_R\) with PBW symbol \(f\). Then \(F-u\)
remains in the intersection and satisfies \(\iota_B(\sigma_d(F-u))=0\).
Since \(\iota_B\) is injective, \(F-u\in F_{d-1}\cA_B(R)\).
Induction, starting with \(F_{-1}=0\), proves the equality.
\end{proof}

\section{Harish-Chandra isomorphism in semisimple rank one}\label{sec:rankone}

\subsection{\texorpdfstring{\(G_R\) invariants of \(\mathcal C_Q(R)\) and \(\cA_Q(R)\)}{G\_R invariants of C\_Q(R) and A\_Q(R)}}

\begin{lemma}\label{lem:levi}
Let \(Q=LU_Q\) be a standard parabolic of \(G\).
For every commutative ring \(R\), evaluation at the identity gives
canonical identifications
\[
\cA_Q(R)^{G_R}\simeq U_R(\fl_R)^{L_R},
\qquad
\mathcal C_Q(R)^{G_R}\simeq\Sym_R(\fl_R)^{L_R},
\]
filtered and graded, respectively. For \(Q=B\), both targets are
\(S_R\). Under these identifications, \(r_{Q,B}\) and \(r^{\rm cl}_{Q,B}\)
become Harish-Chandra projection \(U_R(\fl_R)^{L_R}\to S_R\)
and restriction \(\Sym_R(\fl_R)^{L_R}\to S_R\), respectively.
The maps \(j_{B,R}\) and \(j^{\rm cl}_{B,R}\) induce the same maps
for \(G\).
\end{lemma}
\begin{proof}
A \(G_R\)-invariant element of \(\cA_Q(R)\) is a constant function
on \(G_R\), whose value lies in \(V_Q(R)^{Q_R}\)
by~\eqref{eq:induced-actions}.
The same argument identifies \(\mathcal C_Q(R)^{G_R}\) with
\(\Sym_R(\fg_R/\fu_{Q,R})^{Q_R}\).

Let \(\Delta\) be the simple roots determined by \(B\), and let
\(I\subseteq\Delta\) be the simple roots of \(L\).
Write \(Z(L)\) for the center of \(L\), and choose a cocharacter
\(\eta:\Gm\to Z(L)\) such that
\[
\langle\alpha,\eta\rangle=0\quad(\alpha\in I),
\qquad
\langle\alpha,\eta\rangle>0\quad(\alpha\in\Delta\setminus I).
\]
PBW gives
\[
V_Q(R)\simeq U_R(\fu_{Q,R}^-)\otimes_R U_R(\fl_R).
\]
Its zero \(\eta\)-weight part is \(U_R(\fl_R)\), and every other
weight is negative. Every \(Q_R\)-invariant element is fixed by
\(\eta\), so~\eqref{eq:weights} gives
\(V_Q(R)^{Q_R}\subseteq U_R(\fl_R)\).

Put \(\fq=\Lie(Q)\). By PBW, the \(Q_R\)-equivariant map
\(U_R(\fq_R)\to V_Q(R)\) has image \(U_R(\fl_R)\) and kernel
\(U_R(\fq_R)\fu_{Q,R}\).
The adjoint action of \(Q_R\) on \(\fq_R/\fu_{Q,R}\simeq\fl_R\)
factors through the adjoint action of \(Q_R/U_{Q,R}\simeq L_R\).
The induced action on \(U_R(\fl_R)\) thus gives
\(V_Q(R)^{Q_R}=U_R(\fl_R)^{L_R}\). The argument for 
\(\Sym_R(\fg_R/\fu_{Q,R})\) is identical.

For \(Q=B\), the part of weight zero is \(S_R\).
On these invariant submodules, PBW identifies \(q_{Q,B}\) with
Harish-Chandra projection for \(L\), and its associated graded map
with the projection \(\Sym_R(\fl_R)\to S_R\).
For every simple root of \(L\), the dot reflections computed in
\(L\) and \(G\) agree by~\eqref{eq:dot}.
\end{proof}

Taking \(Q=G\) in Lemma~\ref{lem:comparisons} and then applying
Lemma~\ref{lem:levi} gives
\begin{equation}\label{eq:invariant-embeddings}
\operatorname{res}_{\fh}:\Sym_R(\fg_R)^{G_R}\into S_R,
\qquad
\HC_R:U_R^{G_R}\into S_R,
\end{equation}
with the second map strict filtered. Here
\(\operatorname{res}_{\fh}\) is induced by projection onto \(\fh_R\)
along the root spaces; its image lies in \(S_R^W\) by normalizer
invariance.

Kac and Weisfeiler \cite[\S9.4(d)]{KW} claim a counterexample to
injectivity of Harish-Chandra projection on \(U_k(\fg_k)^{G_k}\)
for \(G_k=\mathrm{SO}_{2n+1,k}\) in characteristic \(2\).
Polo and Tilouine \cite[\S2.6.1, Lemma and Remark]{PT} correct
this claim by proving injectivity and pointing out that the proposed
kernel element is not \(G_k\)-invariant.

\subsection{Concrete calculations in semisimple rank one}

Let \(D\) be a commutative ring, \(c\in D\), and \(u\) an indeterminate
over \(D\). Put \(q=u(u+c)\). Then
\(D[u]=D[q]\oplus uD[q]\). The \(D\)-algebra involution
\(\tau(u)=-u-c\) fixes \(q\),
and the direct sum decomposition gives
\begin{equation}\label{eq:quadratic-involution}
D[u]^\tau
=D[q]\oplus u\bigl(\Ann_D(2)\cap\Ann_D(c)\bigr)[q].
\end{equation}
Indeed, \(a(q)+ub(q)\) is fixed exactly when \(2b=cb=0\).

\begin{proposition}\label{prop:rankone}
Let \(L\) be split connected reductive of semisimple rank one,
with maximal torus \(T\) and positive root \(\alpha\).
For every commutative ring \(R\), Harish-Chandra projection is a
strict filtered isomorphism
\[
\HC_R:U_R(\fl_R)^{L_R}\iso
\Sym_R(\fh_R)^{\langle s_\alpha\rangle,\bullet}.
\]
Moreover, the ordinary restriction map
\[
\operatorname{res}_{\fh}:\Sym_R(\fl_R)^{L_R}
\longrightarrow \Sym_R(\fh_R)^{\langle s_\alpha\rangle}
\]
is a graded isomorphism if and only if \(R[2]=0\) or \(h_\alpha\)
is primitive in \(Y\).
\end{proposition}
\begin{proof}
Put \(K_0=\ker(\alpha:Y\to\Z)\) and
\(A=\Sym_R(K_0\otimes R)\). This subalgebra is fixed by \(L_R\),
and its image in \(U_R(\fl_R)\) is central.
Take the root vectors \(e=e_\alpha\) and \(f=e_{-\alpha}\), with
\([e,f]=h_\alpha\).
Since \(\alpha(h_\alpha)=2\), the image \(\alpha(Y)\) is either
\(\Z\) or \(2\Z\). We begin with the first case.

\smallskip
\noindent\emph{The case \(\alpha(Y)=\Z\).}
Choose \(H\in Y\) with \(\alpha(H)=1\). Then
\[
Y=\Z H\oplus K_0,\qquad h_\alpha=2H+Z
\quad\text{for some }Z\in K_0.
\]
Consequently \(\Sym_R(\fh_R)=A[H]\), with \(H\) an indeterminate
over \(A\). The Lie brackets are
\[
[H,e]=e,\qquad [H,f]=-f,\qquad [e,f]=2H+Z.
\]
Define
\[
q_0=ef+H(H+Z)\in\Sym_R(\fl_R),\qquad
C=fe+H(H+Z+1)\in U_R(\fl_R).
\]
Both \(q_0\) and \(C\) are invariant under the adjoint action of \(L_R\).
To check this integrally, an element \(x_\alpha(t)\) of the root subgroup
\(U_{\alpha,R}\) acts by
\[
H\longmapsto H-te,\qquad e\longmapsto e,\qquad
f\longmapsto f+t(2H+Z)-t^2e.
\]
Substitution fixes \(q_0\). For \(C\), the quadratic terms cancel
and the linear coefficient is
\[
(2H+Z)e-He-eH-eZ-e=[H,e]-e=0.
\]
The same calculation for \(x_{-\alpha}(t)\in U_{-\alpha,R}\), using
\[
H\longmapsto H+tf,\qquad f\longmapsto f,\qquad
e\longmapsto e-t(2H+Z)-t^2f,
\]
shows that it fixes \(q_0\) and \(C\); for \(C\), the linear coefficient
is \([H,f]+f=0\). The torus \(T_R\) fixes \(H,Z\) and acts on \(e,f\)
with opposite weights \(\alpha,-\alpha\), so it fixes both elements.
Since the big cell \(U_{-\alpha}TU_\alpha\) is dense
in \(L\), these identities give invariance over \(\Z\), and hence
over every \(R\).

Localizing the symmetric algebra at \(e\) gives
\[
\Sym_R(\fl_R)[e^{-1}]=A[e^{\pm1},H,q_0].
\]
The adjoint action of \(x_\alpha(t)\) fixes \(A\) pointwise and fixes
\(e,q_0\), while sending \(H\) to \(H-te\).
For an invariant \(F\) in this localization, comparison
of the highest power of \(t\) in \(F(H-te)=F(H)\) shows that \(F\)
is independent of \(H\). Torus invariance then forces the coefficients
of all nonzero powers of \(e\) to vanish.
Since localization at the polynomial variable \(e\) is
injective, we obtain
\begin{equation}\label{eq:rankone-first-image}
\Sym_R(\fl_R)^{L_R}=A[q_0],\qquad
U_R(\fl_R)^{L_R}=A[C].
\end{equation}
For the second equality, a homogeneous invariant leading symbol is
a polynomial in \(q_0\) with homogeneous coefficients in \(A\).
Replacing \(q_0\) by \(C\) lifts it in the same PBW degree;
subtraction and induction prove the assertion.

The dot reflection fixes \(A\) pointwise and sends \(H\) to \(-H-Z-1\).
Multiplication by \(1+Z\) is injective on \(A\): since \(Z\) has degree one,
the lowest nonzero homogeneous component of any \(b\in A\) remains unchanged
in \((1+Z)b\).
Thus~\eqref{eq:quadratic-involution}, with \(D=A\), \(u=H\), and
\(c=Z+1\), gives
\[
A[H]^{s_\alpha,\bullet}=A[H(H+Z+1)].
\]
Since \(\HC_R(C)=H(H+Z+1)\), this proves the Harish-Chandra assertion
in this case.

For ordinary restriction, \(s_\alpha H=-H-Z\), and
\eqref{eq:quadratic-involution}, with \(D=A\), \(u=H\), and \(c=Z\),
gives
\[
A[H]^{s_\alpha}=A[H(H+Z)]\oplus
H\bigl(\Ann_A(2)\cap\Ann_A(Z)\bigr)[H(H+Z)].
\]
The image of restriction is the first summand, since
\(\operatorname{res}_{\fh}(q_0)=H(H+Z)\).
The second summand vanishes if \(R[2]=0\).
It also vanishes when \(h_\alpha\) is primitive. Indeed, \(h_\alpha\)
then belongs to a \(\Z\)-basis of \(Y\), so multiplication by
\(h_\alpha\) is injective on \(\Sym_R(\fh_R)=A[H]\).
For \(b\in\Ann_A(2)\cap\Ann_A(Z)\), we have
\[
h_\alpha b=(2H+Z)b=0,
\]
and hence \(b=0\).
Conversely, suppose that \(h_\alpha\) is nonprimitive and \(R[2]\ne0\).
Since \(\alpha(h_\alpha)=2\), we have \(h_\alpha\in2Y\), hence
\(Z\in2K_0\). Writing \(Z=2Z_0\) shows that \(Z\) annihilates
\(\Ann_A(2)\), so
\[
\Ann_A(2)\cap\Ann_A(Z)=\Ann_A(2)\ne0,
\]
where nonvanishing follows from \(0\ne R[2]\subseteq\Ann_A(2)\).
Thus the second summand is nonzero, and restriction is not surjective.
This proves the Chevalley criterion in the first case.

\smallskip
\noindent\emph{The case \(\alpha(Y)=2\Z\).}
Write \(h=h_\alpha\). For every \(v\in Y\),
\[
v=\frac{\alpha(v)}2h+
\left(v-\frac{\alpha(v)}2h\right),\qquad
v-\frac{\alpha(v)}2h\in K_0.
\]
Thus \(Y=\Z h\oplus K_0\), and the root datum gives
\(L\simeq\mathrm{SL}_2\times T_0\), with cocharacter lattice \(K_0\)
for \(T_0\) \cite[Thm.~6.1.17]{Conrad}. In particular, \(h\) is primitive.

Use the standard embedding
\(L\simeq\mathrm{SL}_2\times T_0\hookrightarrow\mathrm{GL}_2\times T_0\).
In the diagonal subalgebra of
\(\mathfrak{gl}_2\subset\Lie(\mathrm{GL}_2\times T_0)\), put
\[
H=\operatorname{diag}(1,0),\qquad Z=-I,\qquad h=2H+Z.
\]
The group embedding induces the Lie algebra inclusion
\[
\fl=\mathfrak{sl}_2\oplus\Lie(T_0)
\hookrightarrow\mathfrak{gl}_2\oplus\Lie(T_0).
\]
This inclusion splits as a map of \(\Z\)-modules: a retraction is given by
\((X,v)\mapsto(X-\operatorname{tr}(X)H,v)\).
It remains split after base change to any \(R\), so the induced maps
on their symmetric algebras and universal enveloping algebras are injective as well.

We now show that every \(\mathrm{SL}_{2,R}\)-invariant element \(b\)
of \(\Sym_R(\fl_R)\) or \(U_R(\fl_R)\) is
\(\mathrm{GL}_{2,R}\)-invariant.
Let \(D\) be an \(R\)-algebra and \(g\in\mathrm{GL}_2(D)\), and put
\[
D'=D[t]/(t^2-\det g).
\]
In \(D'\), the element \(t\) is a unit with inverse
\(t(\det g)^{-1}\), and \(s=t^{-1}g\) has determinant \(1\).
Thus \(g=(tI)s\) in \(\mathrm{GL}_2(D')\).
By \(\mathrm{SL}_{2,R}\)-invariance, \(s\) fixes \(b_{D'}\),
and the scalar matrix \(tI\) acts trivially by conjugation.
Hence \(g\) fixes \(b_{D'}\).
Since \(D'\) is free over \(D\) with basis \(1,t\), the natural map
\(N\to N\otimes_D D'\) is injective for every \(D\)-module \(N\).
It follows that \(g\) already fixes \(b_D\).
As \(D\) and \(g\) were arbitrary, \(b\) is
\(\mathrm{GL}_{2,R}\)-invariant.

PBW shows that Harish-Chandra projection is compatible with these
embeddings. By the first case, the Harish-Chandra image in
\(A[H,Z]\) of any \(L_R\)-invariant element is fixed by the dot
reflection. This reflection preserves \(A[h]\), since
\[
s_\alpha\bullet H=-H-Z-1,\qquad s_\alpha\bullet Z=Z,\qquad
s_\alpha\bullet h=-h-2.
\]
Putting \(y=h+1\), we obtain
\begin{equation}\label{eq:sl2-dot-invariants}
\HC_R(U_R(\fl_R)^{L_R})
\subseteq A[y]^{y\mapsto-y}
=A[y^2]\oplus y\,\Ann_A(2)[y^2].
\end{equation}

Using \(C,q_0\) from the \(\mathrm{GL}_2\) case, put
\[
C_{\mathrm{SL}_2}:=y^2+4fe=4C+(Z+1)^2,\qquad
Q_{\mathrm{SL}_2}:=h^2+4ef=4q_0+Z^2.
\]
For \(a\in\Ann_A(2)\), one also has \(ay=a(Z+1)\) and \(ah=aZ\).
Since \(C,q_0,Z\) and every element of \(A\) are fixed by
\(\mathrm{GL}_{2,R}\), these identities give
\[
C_{\mathrm{SL}_2},ay\in U_R(\fl_R)^{L_R},\qquad
Q_{\mathrm{SL}_2},ah\in\Sym_R(\fl_R)^{L_R}.
\]
\[
\HC_R(C_{\mathrm{SL}_2})=y^2,\qquad \HC_R(ay)=ay.
\]
Since \(\HC_R\) fixes \(A\) and is multiplicative on invariants,
these identities give equality in~\eqref{eq:sl2-dot-invariants}.
Injectivity and strictness follow
from~\eqref{eq:invariant-embeddings}.
Likewise,
\[
\operatorname{res}_{\fh}(Q_{\mathrm{SL}_2})=h^2,\qquad
\operatorname{res}_{\fh}(ah)=ah,
\]
so ordinary restriction is surjective onto
\[
A[h]^{h\mapsto-h}
=A[h^2]\oplus h\,\Ann_A(2)[h^2].
\]
Its injectivity follows from~\eqref{eq:invariant-embeddings},
completing the proof.
\end{proof}

\subsection{Groups with primitive coroots}

\begin{corollary}\label{cor:primitive}
Theorem~\ref{thm:main} holds when every coroot is primitive.
\end{corollary}
\begin{proof}
Taking \(G_R\)-invariants in Theorem~\ref{thm:parabolic-reconstruction}
and applying Lemma~\ref{lem:levi} and Proposition~\ref{prop:rankone} gives
\[
\HC_R(U_R^{G_R})=\bigcap_i S_R^{s_i,\bullet}=S_R^{W,\bullet}.
\]
The strict embedding~\eqref{eq:invariant-embeddings} and
Lemma~\ref{lem:multiplicativity} give the strict filtered algebra
isomorphism. Naturality is obvious.
\end{proof}

\section{Central extension and completion of the proof}\label{sec:completion}

Corollary~\ref{cor:primitive} proves the Harish-Chandra isomorphism
when every coroot is primitive. We remove this assumption by constructing
a central extension with primitive coroots. Its induced map on symmetric
algebras allows us to lift every dot-invariant polynomial without
increasing its degree. We then prove the criterion for ordinary
Chevalley restriction.

\subsection{A central extension with primitive coroots}

\begin{proposition}\label{prop:cover}
There is a split connected reductive group \(\widetilde G\) whose coroots
are all primitive, together with a smooth faithfully flat homomorphism
\[
\pi:\widetilde G\longrightarrow G,\qquad \ker\pi\simeq\Gm,
\]
whose kernel is central.
Choose a split maximal torus \(\widetilde T\) mapping onto \(T\), put
\(\widetilde{\fh}=\Lie(\widetilde T)\), and identify the two Weyl groups
via \(\pi\). For every commutative ring \(R\), the differential of \(\pi\)
induces a strict filtered surjection
\[
\Sym_R(\widetilde{\fh}_R)^{W,\bullet}
\twoheadrightarrow \Sym_R(\fh_R)^{W,\bullet}.
\]
\end{proposition}
\begin{proof}
We first construct the root datum of \(\widetilde G\).
Let
\[
\widetilde Y=Y\oplus\Z Z,\qquad
\widetilde X=X\oplus\Z\zeta
\]
be dual lattices, with \(\zeta(Z)=1\), \(\zeta(Y)=0\), and \(X(Z)=0\).
For each root \(\alpha\in\Phi\), set \(\widetilde\alpha=(\alpha,0)\)
and define its coroot by
\begin{equation}\label{eq:rho-coroots}
\widetilde h_\alpha=h_\alpha+\rho(h_\alpha)Z.
\end{equation}
Here \(\rho\) is the half-sum of the positive roots. Since
\(\rho(h_\alpha)=1\) for every simple root, \(\rho\) takes integral
values on the coroot lattice, so \(\widetilde h_\alpha\in\widetilde Y\).
The associated reflection \(\widetilde s_\alpha\) acts on \(\widetilde Y\) by
\[
\widetilde s_\alpha(y+aZ)
=s_\alpha y+\bigl(a-\alpha(y)\rho(h_\alpha)\bigr)Z.
\]
This reflection and its dual preserve the integral lattices.
To check the remaining axioms of a root datum, work over \(\Q\) and use
the change of coordinates
\[
y+aZ\longmapsto y+\bigl(a-\rho(y)\bigr)Z.
\]
It sends \(\widetilde h_\alpha\) to \(h_\alpha\) and leaves every
root \(\widetilde\alpha\) unchanged. Thus the axioms reduce to those
for \(G\), and the new Weyl group identifies with \(W\).
For a simple root, the coefficient of \(Z\) in
\(\widetilde h_\alpha=h_\alpha+Z\) is \(1\), so this coroot is primitive.
Every coroot is Weyl conjugate to a simple coroot; hence all coroots
of the new root datum are primitive.

The existence theorem gives a split reductive group
\(\widetilde G\) with this root datum \cite[Thm.~6.1.16(2)]{Conrad}.
The cocharacter \(Z\) is primitive and is annihilated by every root,
so its image is a central subgroup isomorphic to \(\Gm\).
The quotient by this subgroup has cocharacter lattice
\(\widetilde Y/\Z Z=Y\) and the root datum of \(G\).
The quotient and isomorphism theorems therefore give \(\pi\)
\cite[Lem.~3.3.1, Cor.~3.3.5 and Thm.~6.1.17]{Conrad}.
This quotient map is a \(\Gm\)-torsor, hence smooth and faithfully flat.
It identifies the root subgroups and induces the projection
\(Y\oplus\Z Z\to Y\) on cocharacter lattices. Consequently
\(d\pi:\Lie(\widetilde G)\to\fg\) is surjective after every base change:
it identifies the root spaces and sends \(H+aZ\) to \(H\) on the
Cartan subalgebra.

We now prove the assertion about dot invariants. On symmetric algebras,
the map induced by \(d\pi\) is
\[
\pi_S:S_R[Z]\longrightarrow S_R,\qquad Z\longmapsto0.
\]
Set \(t=Z+1\), so that \(S_R[Z]=S_R[t]\) and \(\pi_S(t)=1\).
For each simple root \(\alpha\), formula~\eqref{eq:dot} gives
\[
\widetilde s_\alpha\bullet H
=s_\alpha H-\alpha(H)t,\qquad
\widetilde s_\alpha\bullet t=t
\quad(H\in Y).
\]
In particular, this action is linear in the generators of \(S_R\)
and in \(t\).

Let \(f\in F_dS_R\cap S_R^{W,\bullet}\), and write
\(f=\sum_{m=0}^d f_m\), where \(f_m\) is homogeneous of degree \(m\).
Define
\begin{equation}\label{eq:rho-homogenization}
\widetilde f=\sum_{m=0}^d t^{\,d-m}f_m
\in\Sym_R(\widetilde{\fh}_R).
\end{equation}
Both \(\widetilde s_\alpha\bullet\widetilde f\) and \(\widetilde f\)
are homogeneous of total degree \(d\) in these generators and \(t\).
At \(t=1\), they become \(s_\alpha\bullet f\) and \(f\), which are equal.
A homogeneous polynomial of total degree \(d\) has a unique expression
\(\sum_{m=0}^d t^{d-m}g_m\), with \(g_m\in\Sym_R^m(\fh_R)\);
its value at \(t=1\) therefore determines it. This proves
\(\widetilde s_\alpha\bullet\widetilde f=\widetilde f\) for every
simple reflection, so \(\widetilde f\) is \(W\)-invariant for the dot action.
Finally, \(\pi_S(\widetilde f)=f\), and substituting \(t=Z+1\)
shows that \(\widetilde f\) has polynomial degree at most \(d\).
Thus \(\pi_S\) is strictly surjective on dot invariants.
\end{proof}

\subsection{Completion of the Harish-Chandra proof}

\begin{proof}[Proof of Theorem~\ref{thm:main}]
By Lemma~\ref{lem:levi} and Proposition~\ref{prop:rankone},
\(\HC_R(U_R^{G_R})\) is fixed by every simple dot reflection and
hence lies in \(S_R^{W,\bullet}\). It remains to prove surjectivity
with the required degree bound.

Choose \(\pi:\widetilde G\to G\) as in Proposition~\ref{prop:cover},
and put \(\widetilde{\fg}=\Lie(\widetilde G)\).
Write \(\pi_U=U_R(d\pi)\) for the induced map on enveloping algebras,
and let \(\pi_S:S_R[Z]\to S_R\) be the induced map on the symmetric
algebras of the Cartan subalgebras, given by \(Z\mapsto0\).
The map \(\pi_U\) is \(\widetilde G_R\)-equivariant, where
\(\widetilde G_R\) acts on \(U_R(\fg_R)\) through \(\pi_R\).
Since \(\pi_R\) is faithfully flat, the identity expressing
\(G_R\)-invariance can be checked after pullback along \(\pi_R\).
Thus \(\pi_U\) sends \(\widetilde G_R\)-invariants to \(G_R\)-invariants.

The differential \(d\pi\) identifies the root spaces and sends \(Z\)
to zero. The PBW decompositions defining Harish-Chandra projection
therefore give a commutative square
\[
\begin{array}{ccc}
U_R(\widetilde{\fg}_R)^{\widetilde G_R}
&\xrightarrow{\ \HC_{\widetilde G}\ }&
\Sym_R(\widetilde{\fh}_R)^{W,\bullet}\\
{\scriptstyle\pi_U}\downarrow && \downarrow{\scriptstyle\pi_S}\\
U_R(\fg_R)^{G_R}
&\xrightarrow{\ \HC_G\ }& S_R^{W,\bullet}.
\end{array}
\]
All four maps preserve the filtrations.

Let \(f\in F_dS_R\cap S_R^{W,\bullet}\).
Proposition~\ref{prop:cover} gives a dot-invariant
\(\widetilde f\in F_d\Sym_R(\widetilde{\fh}_R)\) with
\(\pi_S(\widetilde f)=f\).
All coroots of \(\widetilde G\) are primitive, so
Corollary~\ref{cor:primitive} gives
\(\widetilde u\in U_R(\widetilde{\fg}_R)^{\widetilde G_R}\)
of PBW degree at most \(d\), such that
\(\HC_{\widetilde G}(\widetilde u)=\widetilde f\).
Then \(u=\pi_U(\widetilde u)\) belongs to \(F_dU_R\cap U_R^{G_R}\),
and the square gives \(\HC_G(u)=f\). This proves strict surjectivity.

Injectivity follows from~\eqref{eq:invariant-embeddings}, and
multiplicativity follows from Lemma~\ref{lem:multiplicativity}.
Hence \(\HC_R\) is a strict filtered algebra isomorphism.
Its construction by PBW projection commutes with every ring
homomorphism \(R\to R'\), which proves naturality.
\end{proof}

\subsection{Chevalley restriction over arbitrary rings}\label{subsec:chevalley-gluing}

We now turn to ordinary Chevalley restriction. The following lemma
isolates the factors in which a coroot can fail to be primitive.

\begin{lemma}\label{lem:split}
There is a decomposition
\[
G\simeq G_0\times\prod_{k=1}^m J_{n_k},
\]
where every coroot of \(G_0\) is primitive and \(J_n\) is split
adjoint of type \(B_n\). The decomposition may be chosen to carry
\(T\) to the product of maximal tori of the factors.
Here \(B_1=A_1\) and \(B_2=C_2\).
\end{lemma}
\begin{proof}
If all coroots are primitive, take \(G_0=G\). Otherwise, choose a
nonprimitive coroot \(h_\alpha\). Since \(\alpha(h_\alpha)=2\),
we have \(h_\alpha/2\in Y\). In particular,
\(\beta(h_\alpha/2)\in\Z\) for every root \(\beta\), so
\(h_\alpha/2\) belongs to the weight lattice of the dual root system.
Applying \cite[Lem.~1.4]{CR} to the dual root system, we find that
the component containing \(h_\alpha\) has type \(C_n\) and that
\(h_\alpha\) is long. Hence the component of \(\Phi\) containing
\(\alpha\) has type \(B_n\), with \(\alpha\) short; the rank-one
case is included by \(B_1=A_1\).

Write the roots of this component as
\[
\pm\varepsilon_i\quad(1\le i\le n),\qquad
\pm\varepsilon_i\pm\varepsilon_j\quad(1\le i<j\le n),
\]
and let \(H_1,\ldots,H_n\) be the basis of the \(\Q\)-linear span
of the coroots of this component dual to
\(\varepsilon_1,\ldots,\varepsilon_n\), so that
\(\varepsilon_i(H_j)=\delta_{ij}\).
The coroot of \(\varepsilon_i\) is \(2H_i\).
The Weyl group is transitive on the short roots, so the Weyl orbit
of \(h_\alpha/2\) contains every \(H_i\). Hence \(H_i\in Y\).

Since \(\varepsilon_i\in X\), the map
\[
\theta:Y\longrightarrow\bigoplus_{i=1}^n\Z H_i,\qquad
\theta(y)=\sum_{i=1}^n\varepsilon_i(y)H_i
\]
is a homomorphism of lattices. It satisfies \(\theta(H_i)=H_i\), and thus
\[
Y=\left(\bigoplus_{i=1}^n\Z H_i\right)\oplus\ker\theta.
\]
The roots of the chosen component vanish on \(\ker\theta\).
All other roots vanish on \(\bigoplus_i\Z H_i\), and their coroots
lie in \(\ker\theta\). Together with the dual decomposition of \(X\),
this gives a decomposition of the root datum.
The factor of type \(B_n\) has character lattice
\(\bigoplus_i\Z\varepsilon_i\), which is its root lattice, so it is of adjoint type.
The isomorphism theorem gives \(G\simeq J_n\times G'\), compatibly
with maximal tori \cite[Thm.~6.1.17]{Conrad}.
Repeating this construction for \(G'\) gives the required decomposition.
\end{proof}

\begin{samepage}
\begin{proposition}[Chevalley restriction]\label{prop:chevalley-criterion}
For every commutative ring \(R\), Chevalley restriction
\[
\operatorname{res}_{\fh}:
\Sym_R(\fg_R)^{G_R}\longrightarrow S_R^W
\]
is a graded algebra isomorphism if and only if
\[
R[2]=0\qquad\text{or every coroot of \(G\) is primitive in \(Y\)}.
\]
\end{proposition}
\end{samepage}
\begin{proof}
The map \(\operatorname{res}_{\fh}\) is a graded algebra homomorphism and is injective
by~\eqref{eq:invariant-embeddings}. We therefore need only determine
when it is surjective.

\smallskip
\noindent\emph{Primitive coroots.}
Suppose first that every coroot of \(G\) is primitive.
Taking \(G_R\)-invariants in Proposition~\ref{prop:classical-reconstruction}
and applying Lemma~\ref{lem:levi} and the ordinary restriction
statement in Proposition~\ref{prop:rankone} gives
\[
\operatorname{res}_{\fh}\bigl(\Sym_R(\fg_R)^{G_R}\bigr)
=\bigcap_i S_R^{s_i}=S_R^W.
\]
\smallskip
\noindent\emph{The case \(R[2]=0\).}
Assume now that \(R\) has no nonzero \(2\)-torsion, and use the
decomposition in Lemma~\ref{lem:split}.
We first treat one factor \(J_n\), with Cartan subalgebra
\(\fh_n=\bigoplus_{i=1}^n\Z H_i\).
Its Weyl group acts on the \(H_i\) by signed permutations.
Let \(e_j\) denote the \(j\)-th elementary symmetric polynomial, and put
\[
p_j=e_j(H_1^2,\ldots,H_n^2)\qquad(1\le j\le n).
\]
Invariance under each sign change forces the exponents of all \(H_i\)
to be even, since \(R[2]=0\). The remaining polynomial is symmetric
in \(H_1^2,\ldots,H_n^2\). The theorem on symmetric polynomials,
valid over any commutative coefficient ring, therefore gives
\[
R[H_1,\ldots,H_n]^{W(B_n)}=R[p_1,\ldots,p_n].
\]

Apply Proposition~\ref{prop:cover} to obtain
\(\pi_n:\widetilde J_n\to J_n\).
The construction identifies the symmetric algebra of the Cartan
subalgebra of \(\widetilde J_{n,R}\) with \(R[H_1,\ldots,H_n,Z]\).
The map induced by \(d\pi_n\) sets \(Z=0\). Put
\[
u_i=H_i+(n-i)Z,\qquad \widetilde q_i=u_i(u_i+Z).
\]
For \(i<n\), let \(s_i\) be the reflection corresponding to
\(\varepsilon_i-\varepsilon_{i+1}\), and let \(s_n\) correspond to
\(\varepsilon_n\). For the ordinary Weyl action,
formula~\eqref{eq:rho-coroots} gives
\[
s_iH_i=H_{i+1}-Z,\qquad s_iH_{i+1}=H_i+Z,\qquad
s_nH_n=-H_n-Z.
\]
Each reflection fixes \(Z\) and the remaining \(H_j\).
It follows that \(s_i\) exchanges \(u_i\) and \(u_{i+1}\), while
\(s_n\) sends \(u_n\) to \(-u_n-Z\).
Thus the ordinary Weyl action permutes the \(\widetilde q_i\), and
\[
\widetilde p_j=e_j(\widetilde q_1,\ldots,\widetilde q_n)
\]
is \(W\)-invariant and satisfies \(\widetilde p_j|_{Z=0}=p_j\).

All coroots of \(\widetilde J_n\) are primitive. The case already
proved therefore gives an invariant
\(\widetilde F_j\in\Sym_R(\Lie(\widetilde J_{n,R}))^{\widetilde J_{n,R}}\)
whose restriction is \(\widetilde p_j\).
Set \(F_j=\Sym_R(d\pi_n)(\widetilde F_j)\).
The map on symmetric algebras is equivariant, and \(\pi_{n,R}\)
is faithfully flat, so \(F_j\) is \(J_{n,R}\)-invariant by the same
invariance argument used in the proof of Theorem~\ref{thm:main}.
Since \(d\pi_n\) identifies the root spaces and sets \(Z=0\),
\[
\operatorname{res}_{\fh_n}(F_j)
=\widetilde p_j|_{Z=0}=p_j.
\]
This proves surjectivity for \(J_n\).

To pass to the product, put
\(S_{0,R}=\Sym_R(\fh_{0,R})\), where \(\fh_0=\Lie(T\cap G_0)\),
and write \(W_0\) for the Weyl group of \(G_0\).
Let \(p_{k,j}\) denote \(p_j\) for the factor \(J_{n_k}\).
The polynomial algebra \(S_{0,R}\) also has no nonzero \(2\)-torsion.
Applying the calculation for signed permutations with coefficients in
\(S_{0,R}\), and then taking \(W_0\)-invariants, gives
\[
S_R^W=S_{0,R}^{W_0}
\bigl[p_{k,j}:1\le k\le m,\ 1\le j\le n_k\bigr].
\]
Indeed, the expression as a polynomial in the \(p_{k,j}\) is unique,
and \(W_0\) acts only on its coefficients.
Each coefficient lifts to a \(G_{0,R}\)-invariant polynomial by the
case of primitive coroots, and each \(p_{k,j}\) lifts by the calculation
for \(J_{n_k}\). Viewing these polynomials in \(\Sym_R(\fg_R)\)
and taking their products proves surjectivity for \(G\).

\smallskip
\noindent\emph{Failure of surjectivity.}
Finally, suppose that \(G\) has a nonprimitive coroot and that
\(R[2]\ne0\). Choose an adjoint \(B_n\) factor from Lemma~\ref{lem:split}
and a nonzero element \(a\in R[2]\). The polynomial
\[
\ell=a(H_1+\cdots+H_n)
\]
is fixed by signed permutations because \(a=-a\).
The Weyl groups of the other factors act trivially on it, so
\(\ell\in S_R^W\).

Suppose that \(F\in\Sym_R(\fg_R)^{G_R}\) satisfies
\(\operatorname{res}_{\fh}(F)=\ell\).
The adjoint action and \(\operatorname{res}_{\fh}\) both preserve the grading, so
the degree-one component \(F_1\in\fg_R\) is \(G_R\)-invariant and
satisfies \(\operatorname{res}_{\fh}(F_1)=\ell\).
In particular, \(F_1\in\fg_R^{T_R}=\fh_R\).
The map \(\operatorname{res}_{\fh}\) is the identity on \(\fh_R\), so \(F_1=\ell\).
However, \(x_{\varepsilon_1}(t)\in U_{\varepsilon_1,R}\) acts on \(\ell\) by
\[
\Ad(x_{\varepsilon_1}(t))\ell
=\ell-at e_{\varepsilon_1}\ne\ell.
\]
Here the inequality holds in \(\fg_R[t]\): the coefficient
\(ae_{\varepsilon_1}\) is nonzero because the root space
\(R e_{\varepsilon_1}\) is free and \(a\ne0\).
Thus \(\ell\) is not fixed by \(U_{\varepsilon_1,R}\), a contradiction.
The map \(\operatorname{res}_{\fh}\) is therefore not surjective in this case, which proves
the criterion.
\end{proof}

\section*{Acknowledgements}

Funded by the Deutsche Forschungsgemeinschaft (DFG, German Research
Foundation) under Germany's Excellence Strategy -- The Berlin Mathematics
Research Center MATH+ (EXC-2046/1, EXC-2046/2, project ID: 390685689).

The author was supported by the State of Berlin's Elsa Neumann Scholarship.

\end{document}